\documentclass[a4paper, 10pt, english, reqno, dvipsnames, table]{amsart}
\usepackage[utf8]{inputenc}   
\usepackage[T1]{fontenc}      

\usepackage{amsmath}
\usepackage{amsthm}
\usepackage{amssymb}
\usepackage{amsbsy}    % producing more bold mathematics symbols
\usepackage{mathtools}  % more math symbols 
\usepackage{mathrsfs}   % provide \mathscr
\usepackage{amsfonts}

\usepackage{tikz}
\usepackage{tikz-cd}
\usepackage{quiver}
\usepackage{lipsum}

\usepackage{geometry}
\theoremstyle{plain}
\newtheorem{theorem}{Theorem}[section]
\newtheorem{proposition}[theorem]{Proposition}
\newtheorem{lemma}[theorem]{Lemma}
\newtheorem{corollary}[theorem]{Corollary}

\newtheorem{definition-theorem}[theorem]{Definition-Theorem}
\newtheorem{definition-proposition}[theorem]{Definition-Proposition}

\theoremstyle{definition}

\theoremstyle{remark}

\usepackage{paralist}   % adapt the enumerate environment
\usepackage[figurewithin=none]{caption}
\usepackage{graphicx}
\usepackage{float}
\usepackage{xstring}
\usepackage{multirow}
\usepackage{array}
\usepackage{longtable}

\usepackage{xcolor}

\makeatletter
\def\l@subsection{\@tocline{2}{0pt}{2pc}{5pc}{}}
\makeatother

\usepackage{url}
\usepackage[unicode, colorlinks, linktocpage=true, backref=page]{hyperref}
\hypersetup{
	linkcolor = MidnightBlue,
	citecolor = ForestGreen,
	urlcolor = YellowOrange,
	bookmarksnumbered = true
}

\title[Small Central \(L\)-Values of Prime Quadratic Twists of Elliptic Curves under GRH]{Small Central \(L\)-Values of Prime Quadratic Twists of Elliptic Curves under the Generalized Riemann Hypothesis}

\author{Zikang Dong}
\address[Zikang Dong]{ School of Mathematical Sciences, Soochow University, Suzhou 215006, P. R. China}
\email{zikangdong@gmail.com}

\author{Haidong Li}
\address[Haidong Li]{Jiangsu Police Institute. No.48 Shifo Sangong, Pukou District, Nanjing, Jiangsu, 210031, China.}
\email{lihaidong@jspi.cn}

\date{\today}

\begin{document}

\begin{abstract}
Let $E/\mathbb{Q}$ be an elliptic curve with Hasse--Weil $L$-function $L(s,E)$, normalized so that its central point is $s=\frac12$. Using a long resonator method, we prove that for every fixed $0<c<\sqrt{2}$ and all sufficiently large $X$, there exists a fundamental discriminant $d$ with $X\le |d|\le 2X$ and $|d|$ prime such that the quadratic twist $E_d$ has root number $+1$ and
\[
0 \le L(1/2,E_d) \le \exp\left(-c\sqrt{\frac{\log X\,\log\log\log X}{\log\log X}}\right),
\]
conditionally on the Generalized Riemann Hypothesis for Dirichlet $L$-functions and for the Hasse--Weil $L$-functions of elliptic curves over $\mathbb{Q}$. As an application, we obtain upper bounds for the orders of Shafarevich--Tate groups of prime quadratic twists of \(X_0(49)\) under GRH.
\end{abstract}

    \maketitle
    \tableofcontents

\section{Introduction}
The study of extreme values of various \(L\)-functions is a significant topic in analytic number theory, and many methods have been developed for it. Soundararajan \cite{soundararajan2008extreme} introduced the resonator method to study extreme values of the Riemann \(\zeta\)-function, and this method has since been generalized to other \(L\)-functions. Hua and Huang \cite{hua2023extreme} used the resonator method to study extreme values of the Hasse--Weil \(L\)-function of rational elliptic curves \(E/\mathbb{Q}\), with an application to large orders of Shafarevich--Tate groups for a family of elliptic curves under the Birch--Swinnerton-Dyer conjecture.

Recently, Bondarenko and Seip \cite{bondarenko2017large} introduced the long resonator method, which improved Soundararajan's result. In this article, we use the long resonator method to obtain, under the Generalized Riemann Hypothesis, a small-value estimate for prime quadratic twists with a gain of \(\sqrt{\log\log\log X}\) in the scale of the exponent over the almost-prime estimate of Hua and Huang, with an application to the orders of Shafarevich--Tate groups for a family of elliptic curves. In particular, based on a result of Coates--Li--Tian--Zhai \cite{coates2015quadratic}, we obtain upper bounds for the orders of Shafarevich--Tate groups in this family under the Generalized Riemann Hypothesis.

\subsection{Notations and Conventions}
Let \(E\) be an elliptic curve over \(\mathbb{Q}\) with conductor \(N_E\). Denote by \(L(s,E)\) the Hasse--Weil \(L\)-function of \(E\), which is normalized so that the central point is \(s=\frac{1}{2}\), that is, 
\[
L(s,E)=\sum_{n\ge 1}\frac{\lambda(n)}{n^s}=L_{\mathrm{arith}}(s+\frac{1}{2},E)=\sum_{n\ge 1}\frac{a(n)}{n^{s+1/2}}, \quad \Re(s)>1
\]
where at a prime of good reduction, 
\[
a(p)=p+1-\#E(\mathbb{F}_p).
\]
By the Hasse bound at good primes and the Euler factors at bad primes, we have \(|\lambda(p)|\le 2\). The local Euler factors give \(|\lambda(p^j)|\le j+1\), and hence by multiplicativity \(|\lambda(n)|\le  \tau(n)\), where \(\tau(n)\) is the divisor function. By the modularity theorem over \(\mathbb{Q}\), \(L(s,E)\) extends to an entire function of \(s\). Under this normalization, for \(p\nmid N_E\), let \(\alpha_p,\beta_p\) be the roots of \[
x^2-\lambda(p)x+1=0
\]
and hence \(|\alpha_p|=|\beta_p|=1\) and  \(\alpha_p+\beta_p=\lambda(p)\).

Let \(d\in \mathbb{Z}\) be a fundamental discriminant with \((d,2N_E)=1\) and let \(\chi_d(n)=\left(\frac{d}{n}\right)\) be the associated primitive quadratic Dirichlet character. Denote by \(E_d\) the quadratic twist of \(E\). The Hasse--Weil \(L\)-function is defined by
\[
L(s,E_d)=\sum_{n\ge 1}\frac{\lambda(n)\chi_d(n)}{n^s},\quad \Re(s)>1.
\]
Since \((d,2N_E)=1\), we have the conductor \(Q_d:=N_{E_d}=N_Ed^2\) and the root number \(w(E_d)=w(E)\chi_d(-N_E)\in \{1,-1\}\). 

Following Hua and Huang, let \(N_0=\mathrm{lcm}(N_E,8) \). Fix a sign \(\varepsilon \in\{1,-1\}\) and  fix a residue class \(a_0 \mod N_0\)  such that \( a_0 \equiv 1  \mod 4\) and \(a_0\) is coprime with \(N_0\). Now let 
\[
\Omega(a_0, \varepsilon): = \left\{ d : \begin{array}{l} d \text{ is a fundamental discriminant}, \quad (d, 2N_E) = 1, \\ d \equiv a_0 \pmod{N_0}, \quad \varepsilon d > 0, \quad w(E_d) = 1 \end{array} \right\}.
\]
We choose \(\varepsilon\) and \(a_0\) such that \(\Omega(a_0, \varepsilon)\) is not empty in the following.

The following Generalized Riemann Hypothesis will be used in the article.
\begin{description}
    \item[(GRH-D)] Every primitive Dirichlet \(L\)-function satisfies GRH, including the conductor-one case.
    \item[(GRH-E)] For every \(d \in \Omega(a_0, \varepsilon)\) with \(|d|\) prime, the function \(L(s, E_d)\) satisfies GRH.
\end{description}

\subsection{Main Theorem and its applications}
Write $\log_j $ for the $j$-fold iterated natural logarithm and set
\[
 V_X=\sqrt{\frac{\log X\,\log_3X}{\log_2X}}.
\]

Notations as above, we introduce the main theorem in this article.

\begin{theorem}[Main Theorem]\label{thm:main}
Assume \textup{(GRH-D)} and \textup{(GRH-E)}. Let \(0<c<\sqrt{2}\) be fixed.   Then there exists
\(
X_0=X_0(E,a_0,\varepsilon,c)
\)
such that, for every \(X\ge X_0\), the following holds: there exists \(d\in\Omega(a_0,\varepsilon)\) such that \(|d|\) is prime, \(X\le |d|\le 2X\), and
\[
0\le L(1/2,E_d)\le
\exp\left(-c V_X\right).
\]
\end{theorem}

Let \(\omega(d)\) be the number of distinct prime divisors of \(|d|\). Since \(|d|\) is a prime in the theorem, \(\omega(d)=1\). Therefore, we obtain the following corollary. 

\begin{corollary}\label{cor:almostprime}
Under the hypotheses of Theorem~\ref{thm:main}, for every fixed integer $W\ge1$ and $0<c<\sqrt2$,
\[
 \min_{\substack{d\in\Omega(a_0,\varepsilon)\\X\le|d|\le2X\\\omega(d)\le W}}
 L(1/2,E_d)\le\exp(-cV_X)
\]
for all sufficiently large $X$.
\end{corollary}

We now use the main theorem to give an application to Shafarevich-Tate groups. Now let \(A=X_0(49)\) be  the modular curve of level \(49\), which has genus \(1\) and is viewed as an elliptic curve with the cusp at infinity as its origin.

\begin{corollary}\label{cor:sha49}
Assume \textup{(GRH-D)} and \textup{(GRH-E)} for \(A\) and \(\Omega(a_0,\varepsilon)\) with \((a_0,\varepsilon)=(45,+1)\).
Let \(0< c<\sqrt 2\) be a fixed number. For all sufficiently large $X$, there is a prime \(q\) with \(X\le q \le 2X\) and  \(q\equiv45\mod{392}\) such that \(\mathrm{rank}_{\mathbb{Z}}A_q(\mathbb{Q})=0\), its Shafarevich--Tate group \(\mathrm{Sha}(A_q/\mathbb{Q})\) is finite, and
\[
 1\le\#\mathrm{Sha}(A_q/\mathbb{Q})
 \ll_A\sqrt{X}e^{-cV_X}.
\]
\end{corollary}

\begin{proof}
The conductor of \(A\) is \(49\) and its root number is \(+1\),  and by definition \(N_0=392\).   If \(q\equiv45\mod{392}\), then  \(q\equiv1\mod4\) and \(q\equiv 3 \mod 7\). Moreover, we have 
\(\left(\frac{-7}{q}\right)=\left(\frac q7\right)=\left(\frac37\right)=-1\) which means \(q\) is  inert in $\mathbb{Q}(\sqrt{-7})$. Since \(q>0\) and \(q\nmid49\), we have \(w(A_q)=w(A)\chi_q(-49)=1\). By Dirichlet's theorem, there are primes in this class, so \(\Omega(45,+1)\) is not empty.   Theorem~1.2 in Coates--Li--Tian--Zhai \cite{coates2015quadratic} proves \(L(1/2,A_q)\neq 0\), \(\mathrm{rank}_{\mathbb{Z}} A_q(\mathbb{Q})=0\), \(\# \mathrm{Sha}(A_q/\mathbb{Q})<\infty\), and the full BSD formula for every twist \(A_q\) in this progression. Theorem~\ref{thm:main} therefore supplies a prime \(q\) with \(0<L(1/2,A_q)\le e^{-cV_X}\).

Write \(\mathrm{Tor}_q=\#A_q(\mathbb{Q})_{\mathrm{tors}}\) as the order of torsion subgroups, \(\mathrm{Tam}_q=\prod_p c_p(A_q)\) as Tamagawa number, and let $\Omega_q$ be the integral of the absolute value of a N\'eron differential over \(A_q(\mathbb{R})\). The order \(\mathrm{Tor}_q\) is uniformly bounded by Mazur's theorem, and \(\mathrm{Tam}_q\ge1\). The period relation in \cite[p.~375]{coates2015quadratic} gives $\Omega_q\asymp_A q^{-1/2}$. Under our normalization, the BSD formula of rank \(0\) shows 
\[
 \#\mathrm{Sha}(A_q/\mathbb{Q})
 =\frac{L(1/2,A_q)\mathrm{Tor}_q^2}{\Omega_q \mathrm{Tam}_q}
 \ll_A\sqrt q\,L(1/2,A_q)
 \ll_A\sqrt X\,e^{-cV_X},
\]
which proves the corollary.
\end{proof}

\subsection{Strategy of the proof}

We divide our proof of the main theorem into two parts. First, we relate \(L(1/2,E_d)\) with a polynomial, which is used in the resonator method. In this part, we only need the Generalized Riemann Hypothesis for \(E_d\). Second, we employ the long resonator method, where the Generalized Riemann Hypothesis for Dirichlet characters will be used, to prove our main results.

\subsection*{Acknowledgment}
Generative AI was used only for language polishing. The authors take full responsibility for the final content.

\section{A logarithmic upper bound}
Following the method of Soundararajan \cite{soundararajan2009moments},
we establish an upper bound for $\log L(1/2,E_d)$
in terms of a Dirichlet polynomial. We define a polynomial \(P_x(d)\) as follows. Fix a positive constant \(\alpha<1\) and  let  \(x=X^\alpha\). Define
\[
 P_x(d)=\sum_{\substack{p\le x\\p\nmid N_0}}
       \frac{\lambda(p)\chi_d(p)}{\sqrt p}\,w_x(p), \quad w_x(p)=p^{-1/\log x}\frac{\log(x/p)}{\log x}.
\]

\begin{proposition}\label{prop:logupper}

Fix \(0<\alpha<1\) and let \(x=X^\alpha\). Suppose \(d\in\Omega(a_0,\varepsilon )\), \(X\le|d|\le2X\), and GRH holds for \(L(s,E_d)\). Then, uniformly in such $d$,
\[
 \log L(1/2,E_d)\le P_x(d)+O_{E,\alpha}(\log_2X).
\]
Here we set \(\log 0=-\infty\).
\end{proposition}

We introduce more notations that will be used in this section. Let \(\mathcal{L}:=\log x\) for \(x\ge e\) and \(\sigma_0=\frac{1}{2}+\frac{1}{\mathcal{L}}\). Under \textup{(GRH-E)}, all non-trivial zeroes of \(L(s,E_d)\) has form of \(\rho=1/2+i\gamma \). Let \(
\ell_d(s)=\frac{L'(s,E_d)}{L(s,E_d)}\), and hence  \(
-\ell_d(s)=\sum_{n\ge2}\frac{\Lambda_d(n)}{n^s},\) for \(\Re(s)>1\).

\begin{lemma}
    Assume \textup{(GRH-E)}. Define 
    \[
    S_x(d)=
    \sum_{2\le n\le x}
    \frac{\Lambda_d(n)}{n^{\sigma_0}\log n}
    \frac{\log(x/n)}{\mathcal L}, \quad 
    Z_x(d)=
    \frac1{\mathcal L}
    \sum_\rho\int_{\sigma_0}^{\infty}
    \frac{x^{\rho-u}}{(\rho-u)^2}\,du.
    \]
    Then 
\[
\log L(\sigma_0,E_d)
=S_x(d)-\frac{\ell_d(\sigma_0)}{\mathcal L}
+Z_x(d)+O(1/\mathcal L),\quad |Z_x(d)|\le
 \frac{e^{-1}}{\mathcal L}F_d(\sigma_0),
\]
where \(F_d(\sigma)=\sum_{\gamma}\frac{\sigma-1/2}{(\sigma-1/2)^2+\gamma^2}\).
\end{lemma}

\begin{proof}
    For \(u>1/2\), let \[I(u):=\frac{1}{2\pi i}\int_{\Re(w)=1} -\ell_d(u+w)\frac{x^w}{w^2}dw.\]
    On one hand, use \(-\ell_d(s)=\sum_{n\ge2}\frac{\Lambda_d(n)}{n^s},\) and Perron formula, one has
    \[
    I(u)=\sum_{n\le x}\frac{\Lambda_d(n)}{n^u}\log(x/n).
    \]
    On the other hand, we can use residue formula to compute \(I(u)\). Let the rectangle be 
    \[
    \mathcal{R}_{R,T}=\{w: -R-u\le \Re(w)\le 1, \quad -T\le \Im(w)\le T \}
    \]
    and let \(T \to +\infty\) and then \(R\to +\infty\), one has
    \[
    I(u)=-\mathcal{L}\ell_d(u)-\ell_d'(u)-\sum_\rho\frac{x^{\rho-u}}{(\rho-u)^2}-\mathcal T_x(u),
    \]
    where \(\mathcal T_x(u)=\sum_{k\ge0}\frac{x^{-u-1/2-k}}{(u+1/2+k)^2}.\) Here, \(-\mathcal{L}\ell_d(u)-\ell_d'(u)\) is the  residue at \(w=0\), \(-\sum_\rho\frac{x^{\rho-u}}{(\rho-u)^2}\) is the sum of residue at non-trivial zeroes and \( \mathcal{T}_x(u)\) is the sum of residue at trivial zeroes given by Gamma function. Therefore,
    \[
    \sum_{n\le x}\frac{\Lambda_d(n)}{n^u}\log(x/n)=-\mathcal{L}\ell_d(u)-\ell_d'(u)-\sum_\rho\frac{x^{\rho-u}}{(\rho-u)^2}-\mathcal T_x(u).
    \]
    Integrating both side from \(\sigma_0\) to \(\infty\) and divide \(\mathcal{L}\), we have 
    \[
    S_x(d)=\log L(\sigma_0,E_d)+\frac{\ell_d(\sigma_0)}{\mathcal{L}}-Z_x(d)-\frac{1}{\mathcal{L}}\int_{\sigma_0}^\infty \mathcal{T}_x(u)du.
    \]
    Since the last term \(\ll \frac{1}{x\mathcal{L^2}}\le \frac{1}{\mathcal{L}^2}\le \frac{1}{\mathcal L}\), we prove the first claim. 

    Since 
    \[
    \int_{\sigma_0}^{\infty} \left|\frac{x^{\rho-u}}{(\rho-u)^2}\right|du \le\frac{\displaystyle\int_{\sigma_0}^{\infty}x^{1/2-u}\,du}{(\sigma_0-1/2)^2+\gamma^2} = \frac{e^{-1}} {\mathcal L\bigl((\sigma_0-1/2)^2+\gamma^2\bigr)},
    \]
    commute the integral and summation in \(Z_x(d)\), we have \(|Z_x(d)|\le \frac{e^{-1}}{\mathcal L}F_d(\sigma_0)\). Hence we proved this lemma.
\end{proof}

\begin{lemma}
    For \(\frac{1}{2}< \sigma \le \frac{3}{2}\), we have 
    \[
    -\ell_d(\sigma)=\frac{1}{2}\log Q_d-F_d(\sigma)+O(1).
    \]
    If \(L(\frac{1}{2},E_d)\neq 0\), then
    \[
    \log L(\frac{1}{2},E_d)\le \log L(\sigma_0,E_d)+\frac{\log Q_d}{2\mathcal{L}}+O(\frac{1}{\mathcal{L}})
    \]
\end{lemma}

\begin{proof}
    We consider the completed \(L\)-function \(\Lambda(s,E_d)\) of \(L(s,E_d)\), which is defined as 
    \[\Lambda(s,E_d)=\left(\frac{Q^{s/2}}{(2\pi)^s}\right) \Gamma(s+\frac{1}{2}) L(s,E_d),\]
    which has functional equation
    \[
    \Lambda(s,E_d)=w(E_d)\Lambda(1-s,E_d)
    \]
    On one hand, take logarithm and derivative, 
    \[
    \frac{\Lambda'(s,E_d)}{\Lambda(s,E_d)}=\frac{1}{2}\log Q_d-\log (2\pi)+\frac{\Gamma'}{\Gamma}(s+\frac{1}{2})+\ell_d(s).
    \]
    On the other hand, by Hadamard factorization theorem,
    \[
    \Lambda(s,E_d)=e^{A+Bs}\prod_{\rho>0} (1-\frac{s}{\rho})e^{s/\rho}.
    \]
    Let \(\rho \) pairs \(1-\rho\), we have 
    \[
    \Lambda(s,E_d)=e^{A'+B's}(s-\frac{1}{2})^m\prod_{\gamma>} (1+\frac{(s-1/2)^2}{\gamma^2}),
    \]
    where \(m\) is the order of \(\Lambda(s,E_d)\) at \(s=\frac{1}{2}\). By functional equation, we have \(B'=0\). Then take logarithm and derivative, we have 
    \[
    \frac{\Lambda'(\sigma ,E_d)}{\Lambda(\sigma ,E_d)}=\frac{m}{\sigma-1/2}+\sum_{\gamma>0}\frac{2(\sigma-1/2)}{(\sigma-1/2)^2+\gamma^2}=F_d(\sigma).
    \]

    By two formulae of \(\frac{\Lambda'(\sigma ,E_d)}{\Lambda(\sigma ,E_d)}\), we have 
    \[
    -\ell_d(\sigma)=\frac{1}{2}\log Q_d-F_d(\sigma)-\log(2\pi)+\frac{\Gamma'}{\Gamma}(\sigma+\frac{1}{2})
    \]
    and for \(1/2<\sigma<3/2\),
    \[
    -\ell_d(\sigma)=\frac{1}{2}\log Q_d-F_d(\sigma)+O(1).
    \]

    Suppose now \(L(1/2,E_d)\neq 0\), then 
    \[
    \begin{aligned}
    \log L(1/2,E_d)
    &=\log L(\sigma_0,E_d)
    +\int_{1/2}^{\sigma_0}-\ell_d(t)\,dt\\
    &=\log L(\sigma_0,E_d)
    +\frac{\sigma_0-1/2}{2}\log Q_d
    -(\sigma_0-1/2)\log(2\pi)\\
    &\quad
    +\log\Gamma\!\left(\sigma_0+\frac12\right)
    -\int_{1/2}^{\sigma_0}F_d(t)\,dt.
    \end{aligned}
    \]
    Since \(x\ge e\), \( \sigma_0-1/2= \frac{1}{\mathcal{L}}=\frac{1}{\log x} \le 1 \). Then \(-\delta\log(2\pi)+\log\Gamma(1+\delta)=O(\delta).\) Moreover, 
    \(F_d(\sigma)\ge 0\), then 
    \[
    \log L(1/2,E_d)\le \log L(\sigma_0,E_d)+\frac{1}{2\mathcal{L}}\log Q_d +O(\frac{1}{\mathcal{L}}).
    \]
\end{proof}

Now we are ready to prove the proposition. 

\begin{proof}[Proof of Proposition~\ref{prop:logupper}] 
    If \(L(1/2,E_d)=0\), then by our convention on \(\log 0=-\infty\), this proposition holds. Now suppose \(L(1/2,E_d)\neq 0\). By the lemmas above,
     \[
    \begin{aligned}
    \log L(1/2,E_d)
    &\le
    S_x(d)-\frac{\ell_d(\sigma_0)}{\mathcal L}
    +Z_x(d)+\frac{\log Q_d}{2\mathcal L}
    +O(1/\mathcal L)\\
    &\le
    S_x(d)+\frac{\log Q}{\mathcal L}
    +\frac{e^{-1}-1}{\mathcal L}F_d(\sigma_0)
    +O(1/\mathcal L)\\
    &\le
    S_x(d)+\frac{\log Q}{\mathcal L}
    +O(1/\mathcal L),
    \end{aligned}
    \]
    We now first consider \(S_x(d)\). Recall \(S_x(d)=\sum_{2\le n\le x}\frac{\Lambda_d(n)}{n^{\sigma_0}\log n}\frac{\log(x/n)}{\mathcal L}\), and \(\Lambda_d(n)\) is supported at prime powers. 
    
    For \(p \nmid N_0\), \(\Lambda_d(p^j)=(\alpha_p^j+\beta_p^j)\chi_d(p^j)\log p,\) with \(|\alpha_p|=|\beta_p|=1.\) Hence, \(|\Lambda_d(p^j)|\le 2\log p\). Therefore, the  \(|\frac{\Lambda_d(n)}{n^{\sigma_0}\log n}\frac{\log(x/n)}{\mathcal L}|\le \frac{2}{jp^{j\sigma_0}}\) when \(n=p^j,p\nmid N_0\).For \(j=1\), the summation over \(p\) and \(p\nmid N_0\) of is just \(P_x(d)\). For \(j=2\), the summation \(\le _{p\le\sqrt x}
    \frac{1}{p}\ll \log \log x\). For \(j\ge 3\), the summation \(\ll 1\).

    For \(p| N_0\), the summation \(\sum_{p^j\le x} \le \sum_{j\ge 1}\frac{C_p}{jp^{j/2}}\) where \(C_p\) depends on \(E\) and \(p\). Since \(p\ge 2\), the summation \(\sum_{p^j\le x} =O_E(1)\). Since bad primes are finite, then \(\sum_{p^j\le x, p|N_0}=O_E(1)\).

    Combine all these, we have
    \[
    S_x(d)=P_x(d)+O_E(\log\log x)=P_x(d)+O_{E,\alpha}(\log\log X).
    \]

    Moreover,  since \(X\le |d|\le 2X\).
    \[
    \frac{\log Q}{\log x}=\frac{\log N_E +2\log |d|}{\alpha \log X}=\frac{2\log X+O(1)}{\alpha \log X}=\frac{2}{\alpha}+O_{E,\alpha}(\frac{1}{\log X})
    \]
    Therefore, 
    \[
    \log L(1/2,E_d)\le P_x(d)+O_{E,\alpha}(\log\log X)+\frac{2}{\alpha}+ O_{E,\alpha}(\frac{1}{\log X})=P_x(d)+O_{E,\alpha}(\log\log X)
    \]
    which completes our proof.
\end{proof}

\section{The resonator method}\label{sec:resonator}
In this section, we use the long resonator method to find a negative value of \(P_x(d)\). We define
\[
\mathcal Q_X=\{q\text{ prime}:X\le q\le 2X,\ \varepsilon q\in\Omega(a_0,\varepsilon)\}.
\]
For a function \(F:\Omega(a_0,\varepsilon)\to\mathbb R\), let
\[
\mathcal A_X(F)=\sum_{q\in\mathcal Q_X}(\log q)F(\varepsilon q),\quad H_X=\frac{X}{\varphi(N_0)}.
\]

\begin{proposition}\label{prop:resonance}

Assume \textup{(GRH-D)}. Fix
\[
0<\theta<\frac12,\quad 0<\gamma<1,\quad 1<A<\frac1\gamma,\quad 0<\alpha<1-2\theta,
\]
and let \(x=X^\alpha\). Then, for all sufficiently large \(X\), the character sum \(R(d)\) defined below satisfies
\[
D_X>0,\quad \frac{N_X}{D_X}=-(2\gamma\sqrt\theta+o(1))V_X,
\]
where \(D_X=\mathcal A_X(R^2)\) and \(N_X=\mathcal A_X(P_xR^2)\).
\end{proposition}

\begin{lemma}\label{lem:primeortho}
    Assume \textup{(GRH-D)}. Let \(u\ge1\) be coprime to \(N_0\), and suppose every prime divisor of \(u\) is less than \(X\). Then
    \[
    \mathcal A_X(\chi_d(u)) =H_X\mathbf 1_{\{u=\square\}} +O_{N_0}\!\left(X^{1/2}\log^2(2N_0uX)\right)
    \]
    uniformly in \(u\).
\end{lemma}

\begin{proof}
    The root number \(w(E_d)=w(E)\chi_d(-N_E)\) is constant for fundamental discriminants of sign \(\varepsilon\) in the class \(a_0\bmod N_0\). Since \(\Omega(a_0,\varepsilon)\ne\emptyset\), this constant is \(1\). For a prime \(q\equiv\varepsilon a_0\bmod N_0\), the integer \(\varepsilon q\) is odd, squarefree, congruent to \(1\) modulo \(4\), and coprime to \(2N_E\). Thus \(\mathcal Q_X\) consists of the primes in \([X,2X]\) in this class.

    Write \(u=u_1v^2\) with \(u_1\) squarefree. Since \((u,N_0)=1\), the integer \(u_1\) is odd. Also \(q\nmid u\) for \(q\in\mathcal Q_X\), so
    \[
    \chi_{\varepsilon q}(u)=\left(\frac{\varepsilon}{u_1}\right)\left(\frac{q}{u_1}\right).
    \]
    By character orthogonality modulo \(N_0\), we have
    \[
    \mathcal A_X(\chi_d(u))=\frac{(\varepsilon/u_1)}{\varphi(N_0)}\sum_{\chi\bmod N_0}\overline{\chi(\varepsilon a_0)}\sum_{\substack{X\le q\le 2X\\q\text{ prime}}}(\log q)\chi(q)\left(\frac{q}{u_1}\right).
    \]
    Let \(\xi=\chi(\,\cdot\,)(\,\cdot\,/u_1)\), a character modulo \(N_0u_1\). Applying \cite[Theorem~13.7, \textup{(13.20)}]{montgomery2007multiplicative} at \(X\) and \(2X\), with an endpoint term \(O(\log X)\) when \(X\) is prime, we obtain for the inner sum
    \[
    X\mathbf 1_{\{\xi\text{ principal}\}}+O\!\left(X^{1/2}\log^2(2N_0u_1X)\right).
    \]
    Since \((u_1,N_0)=1\), the two factors are the components of \(\xi\) modulo \(N_0\) and \(u_1\). The Jacobi character modulo the odd squarefree integer \(u_1\) is nonprincipal for \(u_1>1\). Thus \(\xi\) is principal precisely when \(u_1=1\) and \(\chi\) is principal. In this case the contribution is \(X/\varphi(N_0)\). Summing the errors and using \(u_1\le u\), we obtain the lemma.
\end{proof}

Let the parameters be as in Proposition~\ref{prop:resonance} and define
\[
K=X^\theta,\quad L_j=\log_j K\ (1\le j\le3),\quad Y=L_1L_2,\quad J=\lfloor L_2^\gamma\rfloor.
\]
Following \cite{bondarenko2017large,darbar2025large}, we define
\[
\mathscr P_k=\{p\text{ prime}:e^kY<p\le e^{k+1}Y\},\quad \mathscr P=\bigcup_{k=1}^J\mathscr P_k,\quad \psi(p)=\sqrt{\frac{L_1L_2}{L_3}}\, \frac{1}{\sqrt p\log(p/Y)}\quad(p\in\mathscr P).
\]
Let \(\mathscr S\) be the set of squarefree positive integers whose prime factors belong to \(\mathscr P\). Define \(r\) multiplicatively on \(\mathscr S\) by \(r(1)=1\) and \(r(p)=-\lambda(p)\psi(p)\). Let
\[
b_p=r(p)^2,\quad M_k=\frac{L_1}{k^2L_3}.
\]
Since \(J=o(L_2)\), we have \(\log p=(1+o(1))L_2\) uniformly for \(p\in\mathscr P\). Hence \(p\le x\) and \(p\nmid N_0\) for all \(p\in\mathscr P\) when \(X\) is sufficiently large.

\begin{lemma}\label{lem:primeweights}
    Uniformly for \(1\le k\le J\), we have
    \[
    \sum_{p\in\mathscr P_k}b_p\le (1+o(1))M_k.
    \]
    Moreover, if
    \[
    C_K=\sum_{p\in\mathscr P} \frac{\lambda(p)^2\psi(p)w_x(p)}{\sqrt p(1+b_p)},
    \]
    then \(C_K=(\gamma+o(1))\sqrt{L_1L_3/L_2}\).
\end{lemma}

\begin{proof}
    By \cite[Lemma~2.3, p.~2759]{hua2023extreme}, for some \(b_E>0\) we have
    \[
    \vartheta_E(t):=\sum_{p\le t}\lambda(p)^2\log p=t+O_E\!\left(te^{-b_E\sqrt{\log t}}\right).
    \]
    Let \(a=e^kY\), \(b=e^{k+1}Y\), \(f(t)=1/(t\log t)\), and \(\Delta(t)=\vartheta_E(t)-t\). By partial summation,
    \[
    \sum_{p\in\mathscr P_k}\frac{\lambda(p)^2}{p}=\int_a^b f(t)\,dt+[\Delta(t)f(t)]_a^b-\int_a^b\Delta(t)f'(t)\,dt.
    \]
    For \(a\le t\le b\), we have \(\log t=(1+o(1))L_2\) uniformly in \(k\), \(tf(t)\ll 1/L_2\), and \(|tf'(t)|\ll f(t)\). Since \(\log(b/a)=1\), the first integral is \((1+o(1))/L_2\), and the two terms involving \(\Delta(t)\) are \(O_E(e^{-b_E\sqrt{L_2}}/L_2)\). Hence
    \[
    \sum_{p\in\mathscr P_k}\frac{\lambda(p)^2}{p}=\frac{1+o(1)}{L_2}
    \]
    uniformly in \(k\). Since \(\log(p/Y)>k\) for \(p\in\mathscr P_k\), we obtain
    \[
    \sum_{p\in\mathscr P_k}b_p\le\frac{L_1L_2}{k^2L_3}\sum_{p\in\mathscr P_k}\frac{\lambda(p)^2}{p}=(1+o(1))M_k.
    \]

    Since \(p>eY\) and \(\log(p/Y)>1\), we have \(b_p\le4/(eL_3)\). Also, \(\log p/\log x=o(1)\) uniformly on \(\mathscr P\), so \(w_x(p)=1+o(1)\). Therefore
    \[
    C_K=(1+o(1))\sqrt{\frac{L_1L_2}{L_3}}\, \sum_{k\le J}\sum_{p\in\mathscr P_k} \frac{\lambda(p)^2}{p\log(p/Y)}.
    \]
    On the \(k\)-th block, \(k<\log(p/Y)\le k+1\), so
    \[
    \frac{1}{k+1}\sum_{p\in\mathscr P_k}\frac{\lambda(p)^2}{p}\le\sum_{p\in\mathscr P_k}\frac{\lambda(p)^2}{p\log(p/Y)}\le\frac{1}{k}\sum_{p\in\mathscr P_k}\frac{\lambda(p)^2}{p}.
    \]
    The sums of \(1/k\) and \(1/(k+1)\) over \(k\le J\) are both \(\gamma L_3+O(1)\). Therefore the double sum is \((\gamma+o(1))L_3/L_2\), which proves the second assertion.
\end{proof}

For \(m\in\mathscr S\), let \(\omega_k(m)=\#\{p\mid m:p\in\mathscr P_k\}\) and \(T_k=AM_k\). Define
\[
\mathscr R=\{m\in\mathscr S:\omega_k(m)<T_k\text{ for }1\le k\le J\}, \quad R(d)=\sum_{m\in\mathscr R}r(m)\chi_d(m).
\]

\begin{lemma}\label{lem:support}
    For all sufficiently large \(X\), we have
    \[
    |\mathscr R|\le K,\quad \max_{m\in\mathscr R}\log m\ll_A \frac{L_1L_2}{L_3}.
    \]
\end{lemma}

\begin{proof}
    Let \(n_k=|\mathscr P_k|\) and \(z_k=e^{-k}/(k^2L_3)\). By the prime number theorem and \(\log(e^kY)\sim L_2\), we have \(n_k\ll e^kL_1\), so \(n_kz_k\ll M_k\). Since \(0<z_k<1\), we have \(z_k^j\ge z_k^{T_k}\) for \(j<T_k\). Thus the number of choices of prime factors from \(\mathscr P_k\) satisfies
    \[
    \sum_{0\le j<T_k}\binom{n_k}{j}\le z_k^{-T_k}(1+z_k)^{n_k}.
    \]
    Since \(m\) is squarefree, its choices of prime factors in the different blocks determine \(m\). Taking the product over \(k\) and using \(\log(1+z_k)\le z_k\), we obtain
    \[
    \begin{aligned}
    \log|\mathscr R| &\le\sum_{k\le J}\bigl(T_k(k+2\log k+\log L_3)+n_kz_k\bigr)\\
    &\le\frac{AL_1}{L_3}\sum_{k\le J}\frac1k +O_A\!\left(\frac{L_1(1+\log L_3)}{L_3}\right) =(A\gamma+o(1))L_1.
    \end{aligned}
    \]
    Since \(A\gamma<1\), we have \(|\mathscr R|\le K\). For \(m\in\mathscr R\), each block contributes less than \(T_k\) prime factors and \(\log p\ll L_2\) on \(\mathscr P\). Hence
    \[
    \log m=\sum_{p\mid m}\log p\ll L_2\sum_{k\le J}T_k\ll_A \frac{L_1L_2}{L_3}.
    \]
\end{proof}

Define
\[
B=\prod_{p\in\mathscr P}(1+b_p),\quad B_{\mathscr R}=\sum_{m\in\mathscr R}r(m)^2,\quad B_p=\sum_{\substack{v,pv\in\mathscr R\\p\nmid v}}r(v)^2 \quad(p\in\mathscr P).
\]

\begin{lemma}\label{lem:truncation}
    We have \(B_{\mathscr R}=(1+o(1))B\) and \(B_p=(1+o(1))B/(1+b_p)\) uniformly for \(p\in\mathscr P\).
\end{lemma}

\begin{proof}
    Expanding \(B\) gives \(\sum_{m\in\mathscr S}r(m)^2=B\). Define a probability measure on \(\mathscr S\) by \(\mathbb P(m)=r(m)^2/B\). Since
    \[
    \frac{r(m)^2}{B}=\prod_{p\mid m}\frac{b_p}{1+b_p}\prod_{\substack{p\in\mathscr P\\p\nmid m}}\frac{1}{1+b_p},
    \]
    the events \(p\mid m\) are independent with probabilities \(\pi_p=b_p/(1+b_p)\). Let \(Z_k=\omega_k(m)\). By Lemma~\ref{lem:primeweights}, \(\mathbb E Z_k=\sum_{p\in\mathscr P_k}\pi_p\le(1+o(1))M_k\) uniformly in \(k\).

    Let \(t=\log A\). By independence,
    \[
    \mathbb E e^{tZ_k}=\prod_{p\in\mathscr P_k}(1+\pi_p(e^t-1))\le \exp\!\left((e^t-1)\mathbb E Z_k\right).
    \]
    Therefore, by Markov's inequality and \(T_k=AM_k\),
    \[
    \begin{aligned}
    \mathbb P(Z_k\ge T_k-1)&\le e^{-t(T_k-1)}\mathbb E e^{tZ_k}\\
    &\le \exp\!\left(-(A\log A-A+1+o(1))M_k+\log A\right).
    \end{aligned}
    \]
    Since \(A\log A-A+1>0\) and \(\min_{k\le J}M_k=L_1/(J^2L_3)\to\infty\), this is \(\ll e^{-c_AM_k}\) for some \(c_A>0\). Summing over \(k\), we obtain
    \[
    \mathbb P(Z_k\ge T_k-1\text{ for some }k)\ll J\exp\!\left(-\frac{c_AL_1}{J^2L_3}\right)=o(1).
    \]
    In particular, \(\mathbb P(m\in\mathscr R)=1-o(1)\), which gives the formula for \(B_{\mathscr R}\).

    Fix \(p\in\mathscr P\) and condition on \(p\nmid m\), an event of probability \((1+b_p)^{-1}=1+o(1)\) uniformly in \(p\). Thus \(Z_k<T_k-1\) for every \(k\) still holds with conditional probability \(1-o(1)\). Multiplication by \(p\) adds one prime factor in its block, so these inequalities imply \(m,pm\in\mathscr R\). Therefore
    \[
    \frac{B_p}{B}=\mathbb P(p\nmid m)\mathbb P(m,pm\in\mathscr R\mid p\nmid m)=\frac{1+o(1)}{1+b_p}
    \]
    uniformly in \(p\).
\end{proof}

\begin{lemma}\label{lem:resonatoraverages}
    Assume \textup{(GRH-D)}. Then
    \[
    \begin{aligned}
    D_X&=H_X B_{\mathscr R}+O(X^{1/2+\theta+o(1)}B),\\
    N_X&=-2H_X\sum_{p\in\mathscr P} \frac{\lambda(p)^2\psi(p)w_x(p)}{\sqrt p}B_p +O(X^{1/2+\theta+\alpha/2+o(1)}B).
    \end{aligned}
    \]
\end{lemma}

\begin{proof}
    Expanding \(R^2\) and \(P_xR^2\), we have
    \[
    \begin{aligned}
    D_X&=\sum_{m,n\in\mathscr R}r(m)r(n)\mathcal A_X(\chi_d(mn)),\\
    N_X&=\sum_{\substack{p\le x\\p\nmid N_0}}\frac{\lambda(p)w_x(p)}{\sqrt p}\sum_{m,n\in\mathscr R}r(m)r(n)\mathcal A_X(\chi_d(pmn)).
    \end{aligned}
    \]
    Since the primes in \(\mathscr P\) are at most \(x<X\) and are coprime to \(N_0\), Lemma~\ref{lem:primeortho} applies to \(mn\) and \(pmn\) in these sums. By Lemma~\ref{lem:support}, \(\log(2N_0pmnX)\ll \log X+L_1L_2/L_3\), so the error in each average is \(X^{1/2+o(1)}\). By Cauchy--Schwarz,
    \[
    \left(\sum_{m\in\mathscr R}|r(m)|\right)^2\le |\mathscr R|B_{\mathscr R}\le KB.
    \]
    Hence the total error in \(D_X\) is \(O(X^{1/2+o(1)}KB)\). For \(N_X\), we also use \(0\le w_x(p)\le1\) and \(|\lambda(p)|\le2\), which give
    \[
    \sum_{\substack{p\le x\\p\nmid N_0}}\frac{|\lambda(p)|w_x(p)}{\sqrt p}\ll\sum_{n\le x}\frac{1}{\sqrt n}\ll\sqrt x.
    \]
    This gives an error \(O(X^{1/2+o(1)}K\sqrt{x}B)\) for \(N_X\).

    Since \(m,n\) are squarefree, \(mn=\square\) exactly when \(m=n\). This gives the main term \(H_X B_{\mathscr R}\) for \(D_X\). For \(pmn=\square\), the prime \(p\) divides exactly one of \(m,n\), and all other prime factors occur in both. Thus \((m,n)=(pv,v)\) or \((v,pv)\), with \(p\nmid v\) and \(p\in\mathscr P\). Since \(r(pv)=-\lambda(p)\psi(p)r(v)\), these two pairs give the main term
    \[
    -2H_X\sum_{p\in\mathscr P}\frac{\lambda(p)^2\psi(p)w_x(p)}{\sqrt p}\sum_{\substack{v,pv\in\mathscr R\\p\nmid v}}r(v)^2.
    \]
    The inner sum is \(B_p\), which completes the proof.
\end{proof}

Now we are ready to prove the proposition.

\begin{proof}[Proof of Proposition~\ref{prop:resonance}]
    Since \(\theta+\alpha/2<1/2\) and \(H_X=X/\varphi(N_0)\), both errors in Lemma~\ref{lem:resonatoraverages} are \(o(H_XB)\). By Lemma~\ref{lem:truncation}, \(B_{\mathscr R}=(1+o(1))B\) and \(B_p=(1+o(1))B/(1+b_p)\) uniformly in \(p\). Since \(\lambda(p)^2\psi(p)w_x(p)\ge0\), it follows that
    \[
    D_X=(1+o(1))H_XB>0,\quad N_X=-2(1+o(1))H_XB C_K+o(H_XB).
    \]
    Since \(K=X^\theta\), we have \(L_1=\theta\log X\), \(L_2=\log_2X+O(1)\), and \(L_3=\log_3X+o(1)\). By Lemma~\ref{lem:primeweights},
    \[
    C_K=(\gamma+o(1))\sqrt{\frac{L_1L_3}{L_2}}=(\gamma\sqrt\theta+o(1))V_X.
    \]
    Since \(V_X\to\infty\), we obtain
    \[
    \frac{N_X}{D_X}=-2(1+o(1))C_K+o(1)=-(2\gamma\sqrt{\theta}+o(1))V_X.
    \]
\end{proof}

\section{The proof of the Main Theorem}\label{sec:mainproof}
Now we are ready to prove the main theorem.

\begin{proof}[Proof of Theorem~\ref{thm:main}]
    Let \(0<c<\sqrt{2}\) be fixed and let \(t=(1+c/\sqrt{2})/2\). Define
    \[
    \theta=\frac{t}{2},\quad \gamma=\sqrt{t},\quad A=\frac{1+1/\gamma}{2},\quad \alpha=\frac{1-t}{2}.
    \]
    Since \(c/\sqrt{2}<t<1\), we have \(0<\theta<\frac12\), \(0<\gamma<1\), and
    \[
    2\gamma\sqrt{\theta}=\sqrt{2}\,t>c,\quad A>1,\quad A\gamma=\frac{1+\gamma}{2}<1,\quad 0<\alpha<1-t=1-2\theta.
    \]
    Therefore Proposition~\ref{prop:resonance} holds for  \(x=X^\alpha\).

    For all sufficiently large \(X\), we have \(D_X>0\). The quotient \(N_X/D_X\) is a weighted average of \(P_x(\varepsilon q)\) with weights \((\log q)R(\varepsilon q)^2\). Since these weights are nonnegative, we obtain
    \[
    \min_{\substack{q\in\mathcal Q_X\\R(\varepsilon q)\ne0}}P_x(\varepsilon q)\le\frac{N_X}{D_X}=-(2\gamma\sqrt{\theta}+o(1))V_X.
    \]
    Let \(q\) attain this minimum and let \(d=\varepsilon q\). Then \(d\in\Omega(a_0,\varepsilon)\), \(|d|\) is prime, and \(X\le |d|\le 2X\).

    By \textup{(GRH-E)}, \(L(s,E_d)\) has no zeros for real \(s>\frac12\). Since \(L(s,E_d)\) is real for real \(s\) and tends to \(1\) as \(s\to+\infty\), it is positive for real \(s>\frac12\). Hence \(L(1/2,E_d)\ge0\) by continuity.

    By Proposition~\ref{prop:logupper} and the bound for \(P_x(d)\), we have
    \[
    \begin{aligned}
    \log L(1/2,E_d)&\le P_x(d)+O_{E,\alpha}(\log_2X)\\
    &\le -(2\gamma\sqrt{\theta}+o(1))V_X+O_{E,\alpha}(\log_2X).
    \end{aligned}
    \]
    Since \(\log_2X/V_X=\sqrt{(\log_2X)^3/(\log X\,\log_3X)}=o(1)\), the last error is \(o(V_X)\). As \(2\gamma\sqrt{\theta}>c\), we have \(\log L(1/2,E_d)\le -cV_X\) for all sufficiently large \(X\). Therefore
    \[
    0\le L(1/2,E_d)\le \exp(-cV_X).
    \]
\end{proof}

\bibliographystyle{alpha}
\bibliography{ref}

\end{document}